\documentclass[preprint,number,12pt]{elsarticle}
\usepackage{amssymb}
\usepackage{amsmath}
\usepackage{subfig}
\usepackage{pgfplots}
\newenvironment{proof}{\noindent{\bf Proof: }}{\hfill\fbox{}\vspace*{3mm}}

\def\b{\bf}

\def\dfrac{\displaystyle \frac}
\def\dsum{\displaystyle \sum}
\def\dint{\displaystyle \int}

\def\b1{{\bf{1}}}
\def\={\!\!\!=\!\!\!}

\begin{document}

\begin{frontmatter}
\title{An inventory model for group-buying auction}
\author{Allen H. Tai}

\address{Department of Applied Mathematics, The Hong Kong Polytechnic University, Hung Hom, Hong Kong}

\begin{abstract}
Group-buying auction has become a popular marketing strategy in the last decade.
In this paper, a stochastic model is developed for an inventory system subjects to demands from group-buying auctions.
The model discussed here takes into the account of the costs of inventory, transportation, dispatching and re-order as well as the penalty cost of non-successful auctions. 
Since a new cycle begins whenever there is a replenishment of products,
the long-run average costs of the model can be obtained by using the renewal theory. 
A closed form solution of the optimal replenishment quantity is also derived.

\end{abstract}


\end{frontmatter}


\section{Introduction}
In the last decade, group-buying auction has become popular in many business, especially with the help of web-based platform.
In this paper, we consider a group-buying auction based on the framework described in \cite{Chen}.
In a group-buying auction, a discounted price is set by the seller and customers can place a bid if they accept the price.
If the number of bidders reaches a predetermined minimum within the cutoff time, then the auction ends and the products will be dispatched to customers.
Otherwise, the auction ends and all bidders will leave the auction.
Through this mechanism, sellers can find opportunities for volume selling such that cost reduction can be realized.
Customers can also get the products at discounted prices.

Research works on group-buying auction started to grow at the beginning of this century.
Van Horn et al. \cite{Van Horn} presented a detailed description of group-buying auction and analyzed its benefits to both suppliers and buyers.
Kauffman and Wang \cite{Kauffman(2002)} examined the business models and the pricing mechanisms used by different Internet-based firms which offer online auctions.
A series of case studies and related analyses were also conducted to compare these business models with other new models for Internet-based selling.
Anand and Aron \cite{Anand} gave a general survey of group-buying practices and then derived the optimal group-buying schedule  under varying conditions of heterogeneity in the demand regimes.
Jing and Xie \cite{Jing} suggested that group-buying benefits the sellers in the way that
customers who placed bids will persuade others to join the buying group.
This helps the sellers to reach new consumers.
There were many theoretical and empirical studies on group-buying suggested that the formation of such buyer groups successfully enhanced the bargaining power of customers, see for example \cite{Dana,Inderst,Kauffman(2010a),Kauffman(2010b),Kauffman(2001),Marvel}.

We remark that most research works on group-buying auction focus on the pricing strategies of sellers or bidding strategies of customers, few studies 
consider the inventory management problems faced by the sellers.
Recently, Chen et al. \cite{Chen(2012)} considered a model which integrates group-buying and regular spot-selling option.
A threshold rationing policy was determined for the selling strategy and three inventory policies are proposed:
the optimal inventory control policy (OIC), the first come, first served policy (FCFS), and
the regular customer blocking policy (RCB).
In this paper, we propose a stochastic model for an inventory system subjects to demands from group-buying auctions.
Stochastic models are commonly used in modeling inventory systems with group dispatching, for an overview see for example \cite{Axsater,Cetinkaya(2005)}.
Ching and Tai \cite{Ching} proposed an optimal integrated replenishment and dispatching policy for a vendor-managed inventory (VMI) system. 
A dispatching decision is made whenever the accumulated load reaches a target size or whenever the time since last dispatch reaches a target time before the target load is consolidated.
Later, Mutlu et al. \cite{Mutlu} proposed an similar analytical model for computing the optimal
time-quantity-based policy for consolidated shipments.
They showed that the optimal time-quantity-based policy outperforms the optimal time-based policy.
Cetinkaya et al. \cite{Cetinkaya} considered the case when vendors are facing demands which arrive randomly in random sizes and
are allowed to consolidate small orders until an economical dispatch quantity accumulates.
Marklund \cite{Marklund} studied the joint effect of inventory replenishment and shipment consolidation.
A supply chain with a central warehouse and $n$ nonidentical retailers were considered and 
an exact recursive method for evaluating the expected cost was also given.

The remaining of the paper is organized as follows. In Section 2, we present the inventory 
model for group-buying auction as well as a cost
analysis of the model and then derive the optimal replenishment quantity.
We give a numerical example in Section 3.
Finally, concluding remarks are given in Section 4 to conclude the paper
and address further research issues.

\section{The inventory model and cost analysis}
In this section, we first describe the mechanism of a group-buying auction and then give a inventory model corresponding to it.
Here we give a general structure of group-buying auction \cite{Chen}.
In a group-buying auction, the seller needs to determine three quantities:
a discounted price for one auctioned product; the maximum auction time for each auction, $T$, and 
the number of bidders required for a successful auction, $N$.
Suppose that when a bidder places a bid, it means he is willing to pay the discounted price set by the seller to buy the  auctioned product.
An auction begins at time $0$ and ends when either of the following two cases happens: 
\begin{enumerate}
	\item[(i)] the number of bidders reach $N$. In this case the auction is success. $N$ auctioned products will be dispatched from the inventory and delivered to the bidders.
	\item[(ii)] time reaches $T$. In this case the auction is not success. All the bids are then canceled.
\end{enumerate}

For simplicity of discussion, we assume that the bidder arrival process is a Poisson process.
Because of the stationary and independent increments properties of the Poisson process,
the expected dispatch time for each dispatch is identical.
The expected number of non-successful auctions before each dispatch is also identical.

The following notation are used in the development of the model.

$T$: the maximum auction time for each auction;

$N$: number of bidders required for a successful auction;

$\lambda$: bidder arrival rate;

$Q$: replenishment quantity;

$D$: dispatching cost;

$F$: the unit transportation cost;

$I$: the unit inventory cost per unit of time;

$C_p$: penalty cost for a non-successful auction;

$K$: the re-order cost.

The followings are the assumptions of the model.
\begin{enumerate}
\item[(A1)] After an auction ends, a new auction starts instantaneously. 
\item[(A2)] The inventory level is under continuously review.
	\item[(A3)] The lead time of inventory replenishment is assumed to be negligible.
	\item[(A4)] At the time of a dispatch, if the inventory in stock is not enough to clear the demand, 
	the replenishment arrives immediately and bring the inventory level back to $Q$.
\end{enumerate}

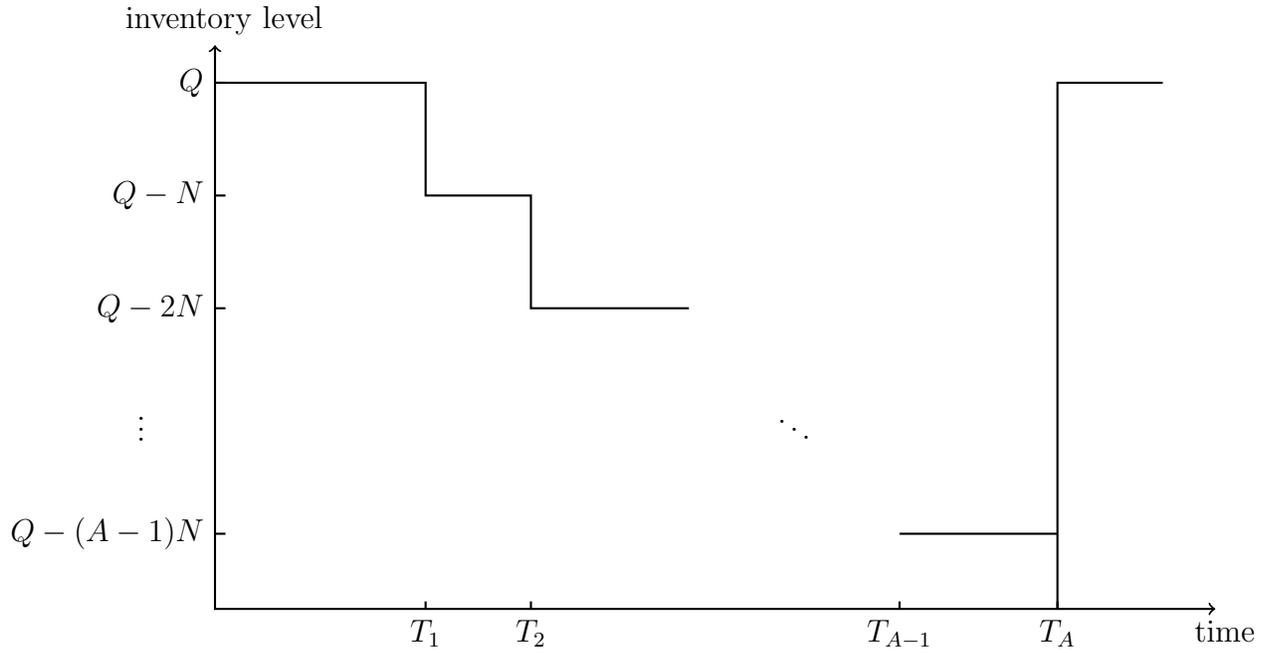
\begin{figure} 
\centering 
\begin{tikzpicture} [ domain=0:10,xscale=1.4,yscale=1] 
\draw [thick,<->] (0,7.5) node [above] {\; inventory level} -- (0,0) -- (0,0) -- (9.5,0) node [below] {\; time};
\draw[thick] (0,7) -- (2,7) -- (2,5.5) -- (3,5.5) -- (3,4) -- (4.5,4);
\draw[thick] (0,7) -- (0.1,7) node [left] {$Q\ $};
\draw[thick] (0,5.5) -- (0.1,5.5) node [left] {$Q-N\ $};
\draw[thick] (0,4) -- (0.1,4) node [left] {$Q-2N\ $};
\node at (5.5,2.5) {$\ddots$};
\node at (-0.7,2.5) {$\vdots$};
\draw[thick] (6.5,1) -- (8,1);
\draw[thick] (8,0) -- (8,7) -- (9,7);
\draw[thick] (0,1) -- (0.1,1) node [left] {$Q-(A-1)N\ $};

\node [below] at (2,0) {$T_1$};
\draw[thick] (2,0) -- (2,0.1);
\node [below] at (3,0) {$T_2$};
\draw[thick] (3,0) -- (3,0.1);
\node [below] at (6.5,0) {$T_{A-1}$};
\draw[thick] (6.5,0) -- (6.5,0.1);
\node [below] at (8,0) {$T_A$};
\draw[thick] (8,0) -- (8,0.1);
\end{tikzpicture} 
\caption{The inventory level of a replenishment cycle.}
\end{figure} 

A realization of the inventory levels in a replenishment cycle is shown in Figure 1.
The aim of this paper is to develop an inventory model in which the demands are from group-buying auctions. 
Based on the model, the optimal values of the replenishment quantity $Q$ can be obtained such that
the average long-run cost of the system is minimized.

\subsection{The inventory model}
Suppose that an auction starts at time $0$.
We let $S_N$ be the time instant that the $N$th bidder arrive.
The probability density function (pdf) of $S_N$ is given by \cite{Ross}
$$
f_{S_N}(t)=\lambda e^{-\lambda t} \dfrac{(\lambda t)^{N-1}}{N!},
$$
and the cumulative distribution function (cdf) is given by
$$
F_{S_N}(t)=1-\sum_{n=0}^{N-1} \dfrac{e^{-\lambda t}(\lambda t)^n}{n!}.
$$

We let $p_T$ be that probability of an auction ends before time $T$, which means there are $N$ bidders before time $T$, then
\begin{equation}\label{pT}
p_T=\Pr(S_N <T)=F_{S_N}(T)=1-\sum_{n=0}^{N-1} \dfrac{e^{-\lambda T}(\lambda T)^n}{n!}.
\end{equation}
Therefore, the expected duration of an auction given that the auction ends before $T$ is
\begin{equation}\label{ETa}
E[T_a]=\dint_0^T t\, \dfrac{ f_{S_N}(t) }{p_T}\, dt = T-\dint_0^T F_{S_N}(t)\, dt=\sum_{n=0}^{N-1} \dfrac{\Gamma(n+1,\lambda T)}{\lambda \cdot n!}.
\end{equation}

Denote by $T_i\ (i = 1, 2, \ldots ,A)$ the instants where a dispatch takes place.
Since the bidder arrival process is a Poisson process, the expected dispatch time for each dispatch is identical.
We let 
$$
T_d=E[T_i-T_{i-1}], \quad \mbox{for }i = 1, 2, \ldots ,A \mbox{ and }T_0=0. 
$$
If there are $N$ bidders before time $T$, then the expected dispatch time is just $E[T_a]$.
Otherwise, the auction starts again.
Hence the expected dispatch time is given by
$$
T_d = p_T E[T_a]+(1-p_T) (T+T_d).
$$
Solving it we get,
\begin{equation}\label{ETd}
T_d=\dfrac{1-p_T}{p_T}T+\sum_{n=0}^{N-1} \dfrac{\Gamma(n+1,\lambda T)}{\lambda \cdot n!}.
\end{equation}
By similar arguments, the expected number of non-successful auctions before a dispatch is given by
$$
E[N_a]=p_T \cdot 0+(1-p_T) (1+E[N_a]).
$$
Solving it we get,
$$
E[N_a]=\dfrac{1-p_T}{p_T}.
$$

If after the $A$th dispatch (for certain $A$) the inventory is out of stock, then a replenishment order is placed and arrived at once as we assume zero lead time.
Therefore, the number of dispatches in a replenishment cycle, $A$, is given by
$$
A=\left\lceil \dfrac{Q}{N}\right\rceil.
$$
Here $\left\lceil x\right\rceil$ is the ceiling function which gives the smallest integer not less than $x$.
For simplicity of the calculations, we use $Q/N$ to approximate $A$.

\subsection{The cost analysis}

We then derive the average long-run cost for the inventory model. 
Since a new inventory cycle begins whenever
there is a replenishment of products,
we can apply the renewal reward theorem \cite{Barlow} to obtain
the average long-run cost:
$$
C(Q)=\dfrac{\mbox{Expected replenishment cycle cost}}{\mbox{Expected length of a replenishment cycle}}.
$$
\begin{enumerate}
	\item[(i)] The expected inventory cost per cycle is given by
$$
\begin{array}{ll}
&I \times \dsum_{i=1}^A \big[ E[T_i-T_{i-1}] \times (Q-(i-1)N) \big]\\
=& I  T_d \Big[ AQ-N \dsum_{i=0}^{A-1} i \Big]\\[5mm]
=& IA  T_d \Big[ Q-\dfrac{(A-1)N}{2} \Big]\\[5mm]
=& IA  T_d \Big( \dfrac{Q+N}{2} \Big).
\end{array}
$$

\item[(ii)] The dispatching cost per cycle is given by $D\times A = DQ/N$.

\item[(iii)] The transportation cost per cycle is given by $F \times A \times N= FQ$.

\item[(iv)] The expected penalty cost per cycle is given by
$$
C_p \times A \times E[N_a] = \dfrac{C_p(1-p_T)Q}{p_T N}.
$$

\item[(v)] The re-order cost per cycle is given by $K$.
\item[(vi)] The expected length of a replenishment cycle is given by
$$
 \dsum_{i=1}^A E[T_i-T_{i-1}] = A \times T_d = \dfrac{Q T_d}{N}.
$$

\end{enumerate}
Hence, the average long-run cost is given by
$$
C(Q)=\dfrac{IQ}{2}+\dfrac{IN}{2}+\dfrac{D}{T_d}+\dfrac{FN}{T_d}+\dfrac{C_p(1-p_T)}{p_T T_d}+\dfrac{KN}{Q T_d}.
$$

The optimal replenishment quantity can be obtained by minimizing the average long-run cost function $C(Q)$. 
We differentiate $C(Q)$ and get
$$
C'(Q) = \dfrac{I}{2}-\dfrac{KN}{Q^2 T_d}.
$$
We note that the cost function $C(Q)$ is strictly convex for positive $Q$ since
$$
C''(Q)=\dfrac{2KN}{Q^3 T_d}>0 \quad \mbox{for }Q>0.
$$
Thus the unique global minimum for positive $Q$ can be obtained by
solving
$$
C'(Q) = \dfrac{I}{2}-\dfrac{KN}{Q^2 T_d}=0.
$$
The optimal replenishment quantity $Q^*$ is then given by
\begin{equation}\label{Qstar}
Q^*=\sqrt{\dfrac{2KN}{IT_d}},
\end{equation}
with the optimal average long-run cost
\begin{equation}\label{CQstar}
C(Q^*)=\sqrt{\dfrac{2IKN}{T_d}}+\dfrac{IN}{2}+\dfrac{D}{T_d}+\dfrac{FN}{T_d}+\dfrac{C_p(1-p_T)}{p_T T_d}.
\end{equation}

\subsection{Remarks}
\begin{enumerate}
	\item The optimal replenishment quantity $Q^*$ is independent of the dispatching cost $D$,
the unit transportation cost $F$ and  the penalty cost for a non-successful auction $C_p$.
	
	\item For the special case $T \to \infty$, by equation (\ref{pT}), $p_T=1$.
Then we have
$$
T_d=E[T_a]=\dint_0^{\infty} t\, f_{S_N}(t)\, dt = \dfrac{N}{\lambda}.
$$
Hence 
$$
Q^*=\sqrt{\dfrac{2KN}{IT_d}}=\sqrt{\dfrac{2K\lambda}{I}},
$$
which is the classical economic order quantity.
\end{enumerate}

\section{Numerical example}
In this section, we give a numerical example to illustrate the model.
Suppose that a group-buying auction is success if there are $100$ bidders, i.e.\ $N=100$.
Each auction will last for $T=7$ days.
We suppose that the arrival rate of bidders is $\lambda=14$ per day.
The costs are given in Table 1.
\begin{table}
\centering
\begin{tabular}{ccccc}
\hline
$D=40$ & $F=4$ & $I=0.02$ & $C_p=10$ & $K=300$\\
\hline
\end{tabular}
\caption{The costs for the auctions.}
\end{table}
By equation (\ref{pT}), we obtain that the probability of an auction ends before time $T$ is $p_T=0.4333$.
Then by equations (\ref{ETa}) and (\ref{ETd}), the expected duration of an auction given that the auction ends before $T$ is $E[T_a]=6.7829$ and 
the expected dispatch time is $T_d=15.9376$.
Hence by equations (\ref{Qstar}) and (\ref{CQstar}), we obtain the optimal replenishment quantity $Q^*=433.8597$ and the optimal average long-run cost $C(Q^*)=37.47802$.
Moreover, if we restrict $Q^*$ to be an integer, we need to consider the integers close to $Q^*$. Consider $\underline{Q^*}=433$ and $\overline{Q^*}=434$.
Substituting them into the average long-run cost function (\ref{CQstar}) yields $C(\underline{Q^*})=37.47804$ and $C(\overline{Q^*})=37.47802$.
Therefore, in order to minimize the average long-run cost, we should adopt the replenishment policy of which brings the inventory level to $434$ units after each replenishment.

We then consider different combinations of $N$, $T$ and repeat the above procedure to obtain the optimal replenishment quantity for these situations.
The corresponding quantities are given in Table 2.
\begin{table}
\centering
\subfloat[]{
\begin{tabular}{lcccc}
\hline
& $T=6$ & $T=7$ & $T=8$\\
\hline
$N=80$ &     0.6834  &  0.9723 &   0.9994\\
$N=100$ &    0.0484  &  0.4333 &   0.8826\\
$N=120$ &    0.0001  &  0.0172 &   0.2368\\
\hline
\end{tabular}}
\qquad
\subfloat[]{
\begin{tabular}{lcccc}
\hline
& $T=6$ & $T=7$ & $T=8$\\
\hline
$N=80$ &        8.3539 & 5.9060 & 5.7192\\
$N=100$ &    123.94 &  15.938 &  8.1614\\
$N=120$ &    46989 &  406.81 &   33.681\\
\hline
\end{tabular}}

\subfloat[]{
\begin{tabular}{lcccc}
\hline
& $T=6$ & $T=7$ & $T=8$\\
\hline
$N=80$ &        536 &  637 &  648\\
$N=100$ &     156 &  434 &  606\\
$N=120$ &       9  &  94 &  327\\
\hline
\end{tabular}}
\qquad
\subfloat[]{
\begin{tabular}{lcccc}
\hline
& $T=6$ & $T=7$ & $T=8$\\
\hline
$N=80$ &  53.9714 &  72.8599 &  74.9544\\
$N=100$ & 9.1673  & 37.4780  & 65.9756\\
$N=120$ & 3.0524  &  5.7391  & 23.8377\\
\hline
\end{tabular}}
\caption{The values of (a) $p_T$, (b) $T_d$, (c) $Q^*$, (d) $C(Q^*)$ for different $N$, $T$.}
\end{table}

\section{Concluding remarks}
In this paper, we consider an inventory system which is subject to demands from group-buying auction.
An analytic inventory model is developed for the replenishment policy.
By using the renewal reward theorem, closed form solution of the optimal replenishment quantity is also obtained.
We consider constant bidders arrival rate in this paper. 
This constraint may be relaxed in the future research.


\end{document}